\documentclass[a4paper]{amsart}

\usepackage{amsthm, amsmath, amssymb, bbm, mathtools}
\usepackage[margin=4cm]{geometry} 
\usepackage[T1]{fontenc}
\usepackage{mathrsfs}
\usepackage{tikz}
\usepackage{tikz-cd}
\usetikzlibrary{backgrounds}
\usepackage{multicol}
\usepackage{xy} \xyoption{all}
\usepackage{color}
\usepackage{enumerate}
\usepackage[colorlinks=true, pdfstartview=FitV, linkcolor=blue, 
citecolor=blue, urlcolor=blue]{hyperref}

\newtheorem*{rmk}{Remark}
\newtheorem*{obs}{Observation}

\DeclareMathOperator{\Mod}{Mod}
\DeclareMathOperator{\Perv}{Perv}
\DeclareMathOperator{\coker}{coker}
\DeclareMathOperator{\Sol}{\mathcal{S}ol}
\DeclareMathOperator{\RHom}{R\mathcal{H}om}
\DeclareMathOperator{\Coh}{Coh}
\DeclareMathOperator{\rh}{\mathrm{rh}}
\DeclareMathOperator{\an}{\mathrm{an}}
\DeclareMathOperator{\Db}{D^b}

\DeclareMathOperator{\R}{R\!}

\newcommand{\D}{\mathcal{D}}
\newcommand{\A}{\mathbb{A}^1}

\title[Topological computation of the Stokes matrices of $\mathbb{P}(1,3)$]{Topological computation of the Stokes matrices\\of the weighted projective line $\mathbb{P}(1,3)$}
\author[A.-L. Sattelberger]{Anna-Laura Sattelberger}
\address[A.-L.~Sattelberger]{Institut f\"ur Mathematik, Universit\"at Augsburg, 86135 Augsburg, Germany}
\email{anna-laura.sattelberger@math.uni-augsburg.de}
\thanks{}

\begin{document}
\begin{abstract}
The localized Fourier--Laplace transform of the Gau{\ss}--Manin system of $f\colon \mathbb{G}_m \to \A,\ x \mapsto x + x^{-3},$ is a $\D_{\mathbb{G}_m}$-module, having a regular singularity at $0$ and an irregular one at $\infty$. By mirror symmetry, it is closely related to the quantum connection of the weighted projective line $\mathbb{P}(1,3)$. Following \cite{DMHS17}, we compute its Stokes multipliers at $\infty$ by purely topological methods. We compare it to the Gram matrix of the Euler--Poincar\'{e} pairing on $\Db(\Coh(\mathbb{P}(1,3)))$.
\end{abstract}

\maketitle

\setcounter{tocdepth}{1}
\tableofcontents
\section*{Introduction}
In \cite{DMHS17}, A.~D'Agnolo, M.~Hien, G.~Morando, and C.~Sabbah describe how to compute the Stokes multipliers of the enhanced Fourier--Sato transform of a perverse sheaf on the affine line by purely topological methods. To a regular singular holonomic\linebreak $\D$-module $\mathcal{M} \in \Mod_{\rh} \left(\D_{\A}\right)$ on the affine line one associates a perverse sheaf via the regular Riemann--Hilbert correspondence $$\RHom_{\D_{\A}^{\mathrm{an}}} \left((\bullet)^{\an}, \mathcal{O}_{\A}^{\mathrm{an}}\right)[1] \colon \Mod_{\rh} \left(\D_{\A}\right) \stackrel{\simeq}{\longrightarrow} \Perv\left(\mathbb{C}_{\A}\right).$$ 
Let $\Sigma \subset \A$ denote the set of singularities of $\mathcal{M}$. Following \cite[Section~4.2]{DMHS17}, after suitably choosing a total order on $\Sigma$, the resulting perverse sheaf $F\in \Perv_{\Sigma}\left(\mathbb{C}_{\A}\right)$ can be described by linear algebra data, namely the quiver $$\left( \Psi (F),\Phi_{\sigma}(F), u_{\sigma},v_{\sigma} \right)_{\sigma \in \Sigma},$$ where $\Psi(F)$ and $\Phi_{\sigma}(F)$ are finite dimensional $\mathbb{C}$-vector spaces and $u_{\sigma} \colon \Psi(F) \to \Phi_{\sigma}(F)$ and $v_{\sigma} \colon \Phi_{\sigma}(F) \to \Psi(F)$ are linear maps such that $1-u_{\sigma}v_{\sigma}$ is invertible for any $\sigma$. The main result in \cite{DMHS17} is a determination of the Stokes multipliers of the enhanced Fourier--Sato transform of $F$ and therefore of the Fourier--Laplace transform of $\mathcal{M}$ in terms of the quiver of $F$. 

Mirror symmetry connects the weighted projective line $\mathbb{P}(1,3)$ with the Landau--Ginzburg model $$\left( \mathbb{G}_m, f=x+\frac{1}{x^3} \right).$$ The quantum connection of $\mathbb{P}(1,3)$ is closely related to the Fourier--Laplace transform of the Gau{\ss}--Manin system $H^0\left(\int_f \mathcal{O}\right)$ of $f$. We compute that $$F\coloneqq \R f_{\ast}\mathbb{C}[1]\in \Perv_{\Sigma}\left(\mathbb{C}_{\A}\right),$$ where $\Sigma$ denotes the set of singular values of $f$, is the perverse sheaf associated to $H^0\left(\int_f \mathcal{O}\right)$ by the Riemann--Hilbert correspondence. In Section \ref{FLGM}, we compute the localized Fourier--Laplace transform of $f$. In Section \ref{StokesTop}, analogous to the examples in \cite[Section~7]{DMHS17}, we carry out the topological computation of the Stokes multipliers of the Fourier--Laplace transform of $H^0\left(\int_f \mathcal{O}\right)$. In Section \ref{Gram}, we compare the Stokes matrix $S_{\beta}$ that we obtained from our topological computations to the Gram matrix of the Euler--Poincar\'{e} pairing on $\Db(\Coh(\mathbb{P}(1,3)))$ with respect to a suitable full exceptional collection. Following Dubrovin's conjecture about the Stokes matrix of the quantum connection, proven for the weighted projective space $\mathbb{P}\left(\omega_0,\ldots,\omega_n\right)$ in \cite{TU13} by S.~Tanab\'{e} and K.~Ueda and in \cite{CMvdP15} by J.~A.~Cruz Morales and M.~van der Put, they are known to be equivalent after appropriate modifications. We give the explicit braid of the braid group $B_4$ that deforms the Gram matrix into the Stokes matrix $S_{\beta}$.

\section{Gau{\ss}--Manin system and its Fourier--Laplace transform}\label{FLGM}
Let $X$ be affine and $f$ a regular function $f \colon X \to \A$ on $X$. Denote by $\int_f (\bullet )$ the direct image in the category of $\D$-modules and by $M\coloneqq H^0 \left(\int_f \mathcal{O}_X\right)\in \Mod_{\rh}\left(\D_{\A}\right)$ the zeroth cohomology of the Gau{\ss}--Manin system of $f$. Following \cite[Section~2.c]{DS03}, it is given by $$ M= \Omega ^n (X)\left[\partial_t \right] / (d-\partial_t df \wedge )\Omega ^{n-1} (X)\left[\partial_t \right].$$
Denote by $G\coloneqq \widehat{M}\left[\tau ^{-1}\right]$ the Fourier--Laplace transform of $M$, localized at $\tau = 0$. It is given by
\begin{align*}
G = \Omega^n(X)\left[\tau,\tau^{-1}\right] / \left(d - \tau df\wedge \right) \Omega^{n-1}(X)\left[\tau,\tau^{-1}\right].
\end{align*}
$G$ is a free $\mathbb{C} [\tau, \tau^{-1}]$-module of finite rank. Rewriting in the variable $\theta=\tau ^{-1}$ gives the $\mathbb{C} \left[\theta,\theta^{-1}\right]$-module
\begin{align*}
G = \Omega^n(X)\left[\theta,\theta^{-1}\right] / \left(\theta d - df\wedge \right) \Omega^{n-1}(X)\left[\theta,\theta^{-1}\right].
\end{align*}
$G$ is endowed with a flat connection given as follows. For $\gamma = \left[\sum_{k \in \mathbb{Z}} \omega_k \theta ^k \right] \in G$, where $\Omega^n(X)\ni \omega_k =0$ for almost all $k$, the connection is given by (cf. \cite[Definition~2.3.1]{GS}):
\begin{align*}
\theta ^2 \nabla _{\! \frac{\partial}{\partial \theta}} \left(\gamma\right)= \left[\sum_k f \omega_k \theta ^k + \sum_k k \omega_k \theta^{k+1} \right] .
\end{align*}
It is known that $(G,\nabla)$ has a regular singularity at $\theta = \infty$ and possibly an irregular singularity at $\theta = 0$. Rewriting in $\tau=\theta^{-1}$ yields the irregular singularity at $\tau = \infty$.
\\We now consider the Laurent polynomial $f=x+x^{-3}\in \mathbb{C}\left[x,x^{-1}\right]$, being a regular function on the multiplicative group $\mathbb{G}_m$. 
For our computations we pass to the variable $\theta=\tau^{-1}$. We compute that for the given $f$, $G$ is given by the free $\mathbb{C} \left[\theta,\theta^{-1}\right]$-module
 \begin{align*}
 G=\mathbb{C} \left[x,x^{-1}\right]dx\left[\theta,\theta^{-1}\right]/\left( \theta d-\left( dx-3x^{-4}dx\right) \wedge \right)\mathbb{C} \left[x,x^{-1}\right]\left[\theta,\theta^{-1}\right]
 \end{align*} 
 with basis over $\mathbb{C} \left[ \theta,\theta^{-1} \right]$ given by $\left[ \frac{dx}{x}\right], \left[\frac{dx}{x^2}\right],\left[\frac{dx}{x^3}\right], \left[ \frac{dx}{x^4}\right]$. In this basis, the connection is given by 
 \begin{align}\label{GM}
 \theta\nabla_{\! \frac{\partial}{\partial \theta} }= \theta \partial_{\theta} + \begin{pmatrix}
0 & \frac{4}{3\theta} & 0 & 0\\
0 & \frac{1}{3} & \frac{4}{3\theta} & 0\\
0 & 0 & \frac{2}{3} & \frac{4}{3\theta}\\
\frac{4}{\theta} & 0 & 0 & 1
\end{pmatrix}.
 \end{align} 
Via the cyclic vector $m=(1,0,0,0)^{\mathrm{t}}$, we compute the relation 
\begin{align*} 
\nabla_{\theta \partial_{\theta}}^4m + 4 \nabla_{\theta \partial_{\theta}}^3m + \frac{32}{9}\nabla_{\theta \partial_{\theta}}^2m-\frac{256}{27\theta^4}m =0 
\end{align*} and therefore associate the differential operator 
\begin{align*}
P = \left(\theta \partial_{\theta}\right)^4 + 4 \left(\theta \partial_{\theta}\right)^3 + \frac{32}{9}\left(\theta \partial_{\theta}\right)^2 - \frac{256}{27\theta^4}\in \mathbb{C}\left [\theta, \theta^{-1} \right]\langle \partial_{\theta} \rangle = \D_{\mathbb{G}_m}.
\end{align*} 
	
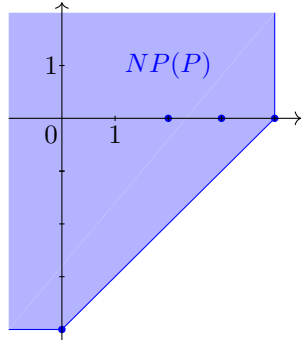
\begin{figure}[h]\label{NP}
\centering
\begin{tikzpicture}[scale=.7]
\filldraw[draw=blue, fill=blue!30] (-1,-4) -- (0,-4) -- (4,0) -- (4,2);
\filldraw[draw=none, fill=blue!30] (-1,2) -- (-1,-4) -- (4,2);
\fill[blue] (0,-4) circle[radius=2pt];
\fill[blue] (2,0) circle[radius=2pt];
\fill[blue] (3,0) circle[radius=2pt];
\fill[blue] (4,0) circle[radius=2pt];
\draw[->] (-1,0) -- (4.5,0);
\draw[->] (0,-4.2) -- (0,2.2);
\draw[] (1,-0.05) -- (1,0.05);
\draw[] (2,-0.05) -- (2,0.05);
\draw[] (3,-0.05) -- (3,0.05);
\draw[] (4,-0.05) -- (4,0.05);
\draw[] (-.05,1) -- (.05,1);
\draw[] (-.05,-1) -- (.05,-1);
\draw[] (-.05,-1) -- (.05,-1);
\draw[] (-.05,-2) -- (.05,-2);
\draw[] (-.05,-3) -- (.05,-3);
\node[blue] at (2,1) {$NP(P)$};
\node[black] at (-.2,-.3) {$0$};
\node[black] at (-.2,1) {$1$};
\node[black] at (1,-.3) {$1$};
\end{tikzpicture}
\caption{Newton polygon of $P$}
\label{NPhyper}
\end{figure}
\noindent The Newton polygon in Figure \ref{NP} confirms that $P$---and therefore system \eqref{GM}---has the nonzero slope $1$ and therefore is irregular singular at $\theta = 0$ and regular singular at $\theta = \infty$.

\section{Topological computation of the Stokes matrices}\label{StokesTop}
We consider the Laurent polynomial $f=x+x^{-3}\colon \mathbb{G}_m \to \A $. Its critical points are given by $\left\{ \pm \sqrt[4]{3}, \pm i \sqrt[4]{3} \right \}$. The critical values of $f$ are given by $$ \Sigma = \left\{ \pm \frac{4}{\sqrt[4]{27}}, \pm \frac{4i}{\sqrt[4]{27}} \right\} \subset \A.$$ The preimages of~\begin{itemize}
	\item $\frac{4}{\sqrt[4]{27}}$ are $ \sqrt[4]{3}$ (double), $\frac{-1-\sqrt{2}i}{\sqrt[4]{27}}$ and $\frac{-1+\sqrt{2}i}{\sqrt[4]{27}}$,
	\item $-\frac{4}{\sqrt[4]{27}}$ are $ - \sqrt[4]{3}$ (double), $\frac{1-\sqrt{2}i}{\sqrt[4]{27}}$ and $\frac{1+\sqrt{2}i}{\sqrt[4]{27}}$,
	\item $i\frac{4}{\sqrt[4]{27}}$ are $ i \sqrt[4]{3}$ (double), $\frac{-i -\sqrt{2}}{\sqrt[4]{27}}$ and $\frac{-i +\sqrt{2}}{\sqrt[4]{27}}$,
	\item $-i\frac{4}{\sqrt[4]{27}}$ are $ -i \sqrt[4]{3}$ (double), $\frac{i +\sqrt{2}}{\sqrt[4]{27}}$ and $\frac{i -\sqrt{2}}{\sqrt[4]{27}}$.
\end{itemize} Since $f$ is proper, we compute by the adjunction formula that 
\begin{align*}
\RHom_{\D^{\an} }\left( \left(\int_f \mathcal{O}\right)^{\an} ,\mathcal{O}^{\an} \right) \simeq \R f^{\an}_{\ast}\RHom_{\D^{\an} } (\mathcal{O}^{\an} ,f^{\dagger}\mathcal{O}^{\an}) \simeq \R f^{\mathrm{an}}_{\ast}\mathbb{C}.
\end{align*}
 Since $f$ is semismall, $\R f_{\ast} \mathbb{C}[1]\in \Perv\left(\mathbb{C}_{\A}\right)$ is a perverse sheaf (cf. \cite{DCM02}). Outside of $\Sigma$, $f$ is a covering of degree $4$, therefore $\R f_{\ast} \mathbb{C}[1]\in \Perv_{\Sigma}\left(\mathbb{C}_{\A}\right)$. By the regular Riemann--Hilbert correspondence 
 \begin{align*}
 \Sol(\bullet)[\dim X] \coloneqq \RHom_{\D_X^{\mathrm{an}}}\left(\left(\bullet\right)^{\an}, \mathcal{O}^{\mathrm{an}}_X\right)\left[\dim X\right] \colon \Mod_{\rh}\left(\D_X\right) \stackrel{\simeq}{\longrightarrow} \Perv (\mathbb{C}_{X^{\an}}),
 \end{align*}
 we associate to $H^0 \left(\int_f \mathcal{O}\right)$ the perverse sheaf $F \coloneqq \R f_{\ast}\mathbb{C}[1]$. 

\begin{figure}[h]
	\includegraphics[width=9.5cm]{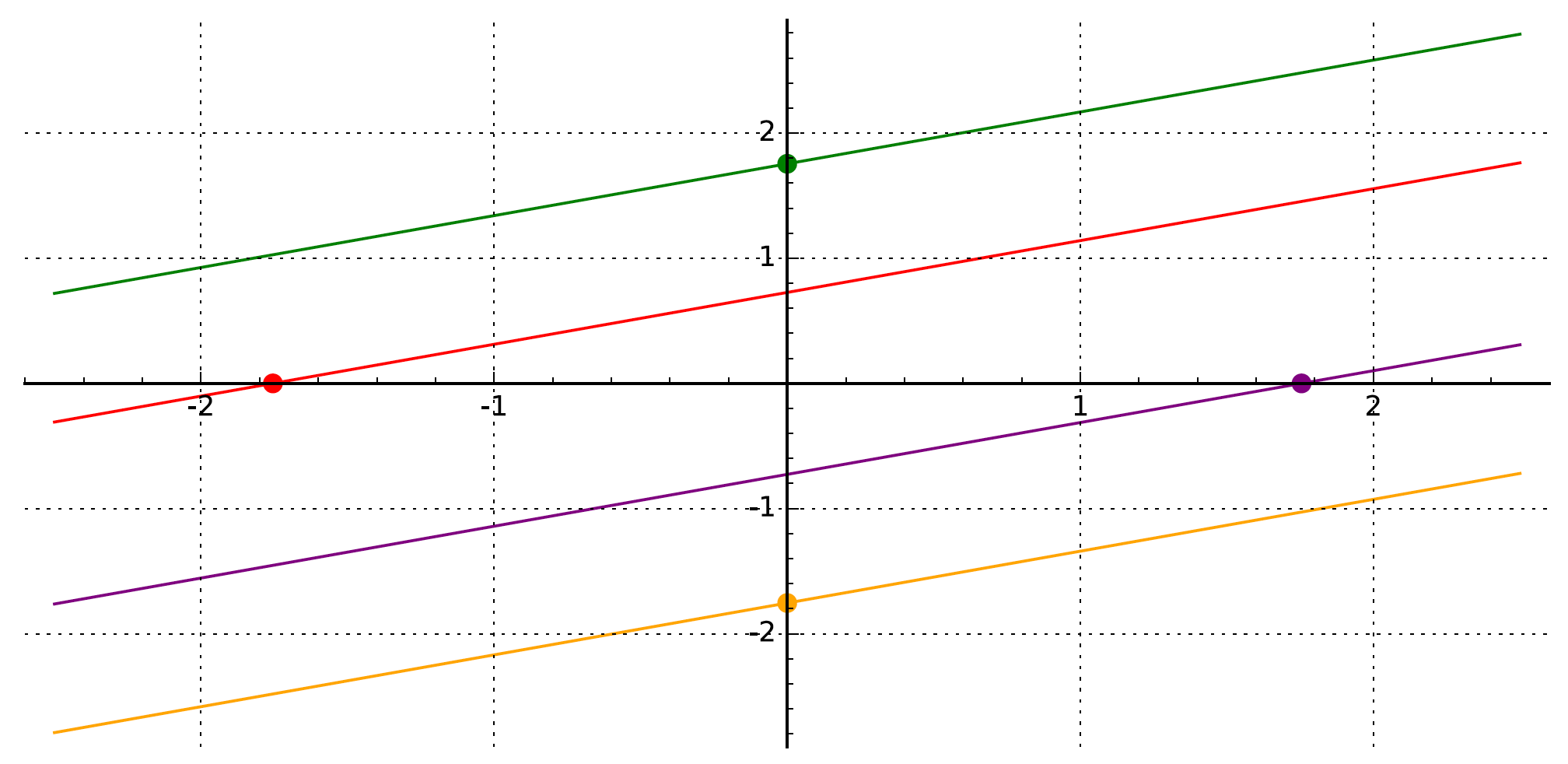}
	\caption{Lines passing through $\sigma_i$ with phase $\frac{\pi}{8}$}
	\label{linesP13}
\end{figure}

\noindent We fix $\alpha = e^{\frac{\pi i}{8}}\in \A,\ \beta = e^{\frac{3 \pi i}{8}}\in \left(\A\right)^{\vee}$, such that $\Re(\langle \alpha, \beta \rangle)=0, \ \Im(\langle \alpha, \beta \rangle)=1$. 
This induces the following order on $\Sigma$ (cf. \cite[Section 4]{DMHS17}): $$\sigma_1\coloneqq \frac{4i}{\sqrt[4]{27}} <_{\beta} \sigma_2\coloneqq -\frac{4}{\sqrt[4]{27}} <_{\beta} \sigma_3\coloneqq \frac{4}{\sqrt[4]{27}} <_{\beta} \sigma_4\coloneqq -\frac{4i}{\sqrt[4]{27}}.$$
In Figure \ref{linesP13}, the $\sigma_i$ are depicted in the following colors:
\begin{center} $\bullet$ \textcolor{green}{$\sigma_1$}: green,\quad $\bullet$ \textcolor{red}{$\sigma_2$}: red,\quad $\bullet$ \textcolor{violet}{$\sigma_3$}: purple, \quad $\bullet$ \textcolor{orange}{$\sigma_4$}: orange.\end{center}

\begin{figure}[h]
	\begin{center}
		\includegraphics[width=6cm]{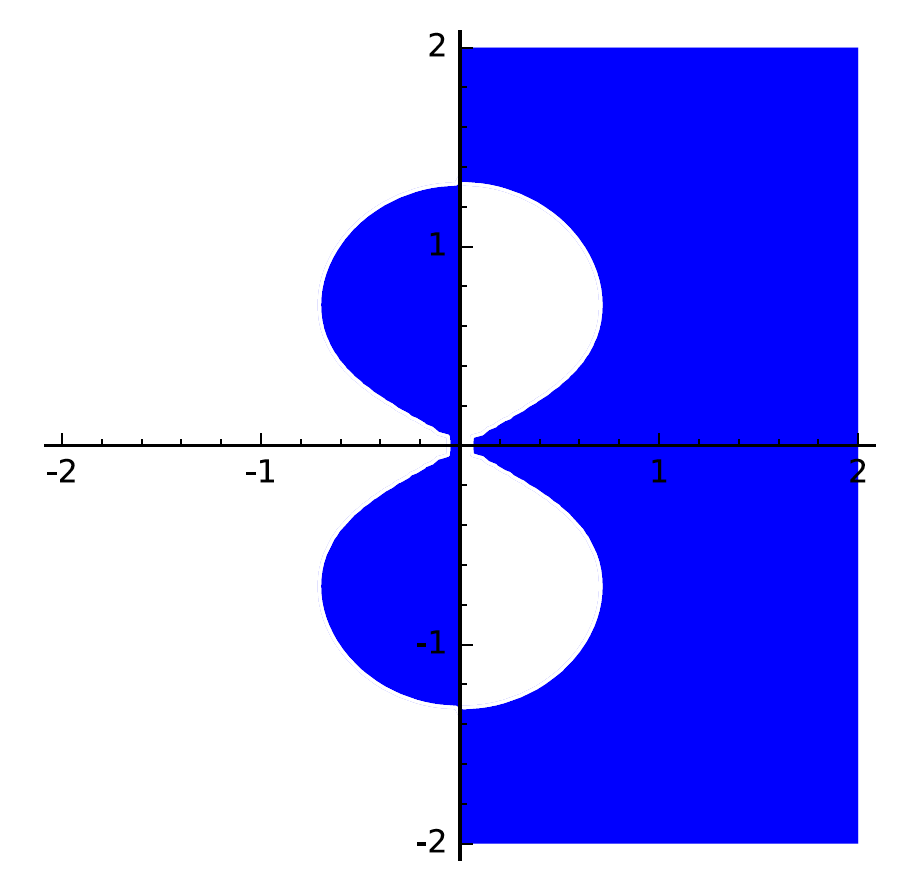}
		\	\includegraphics[width=5.75cm]{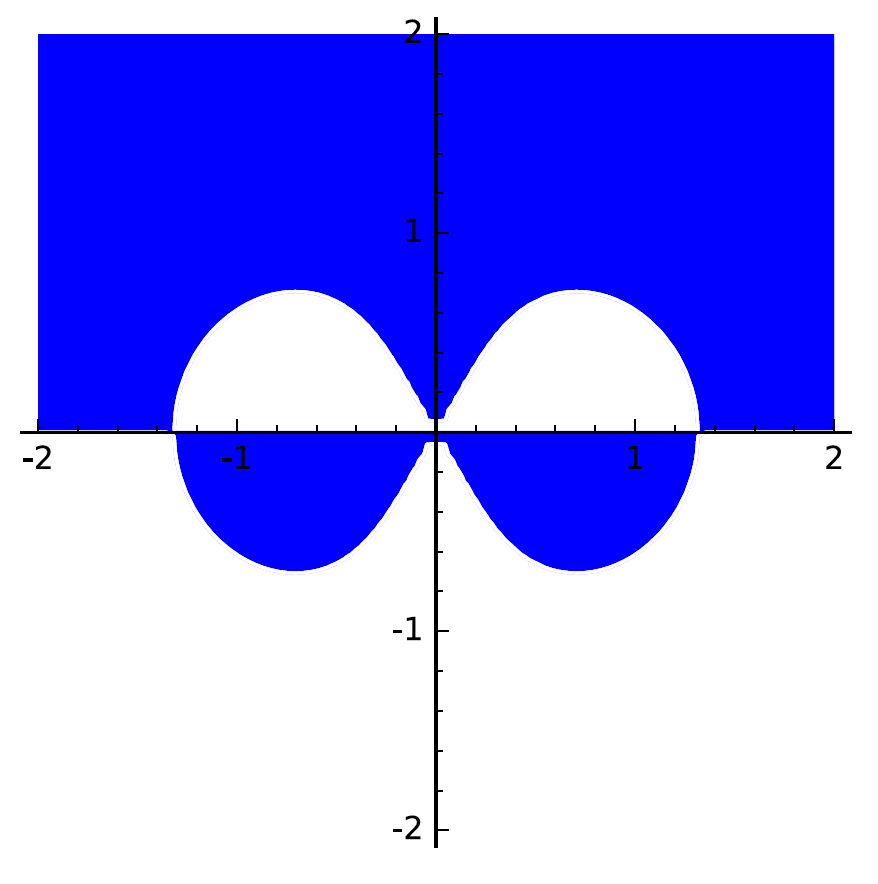}
		\caption{LHS: $\{ x \mid \Re (f(x))\geq 0 \}$, RHS: $\{ x \mid \Im (f(x))\geq 0 \}$}
		\label{realimaginarypartyP13}
	\end{center}
\end{figure}
\begin{figure}[h]
	\includegraphics[width=6cm]{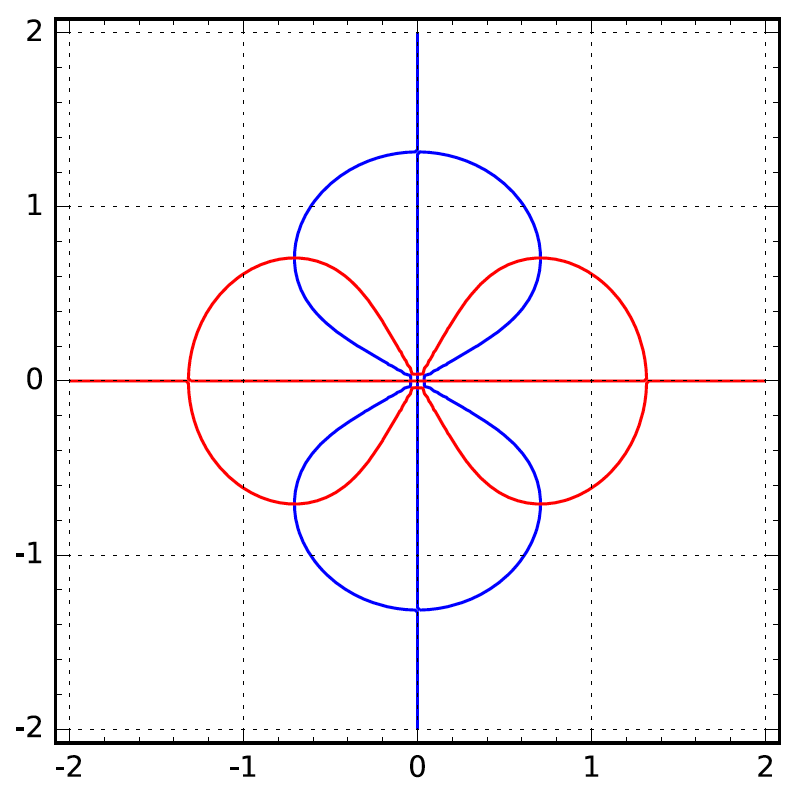}
	\caption{Preimage of the imaginary (resp. real) axis in blue (resp. red) color}
	\label{preimagesaxesP13}
\end{figure}
\begin{figure}[h]
	\includegraphics[width=8cm]{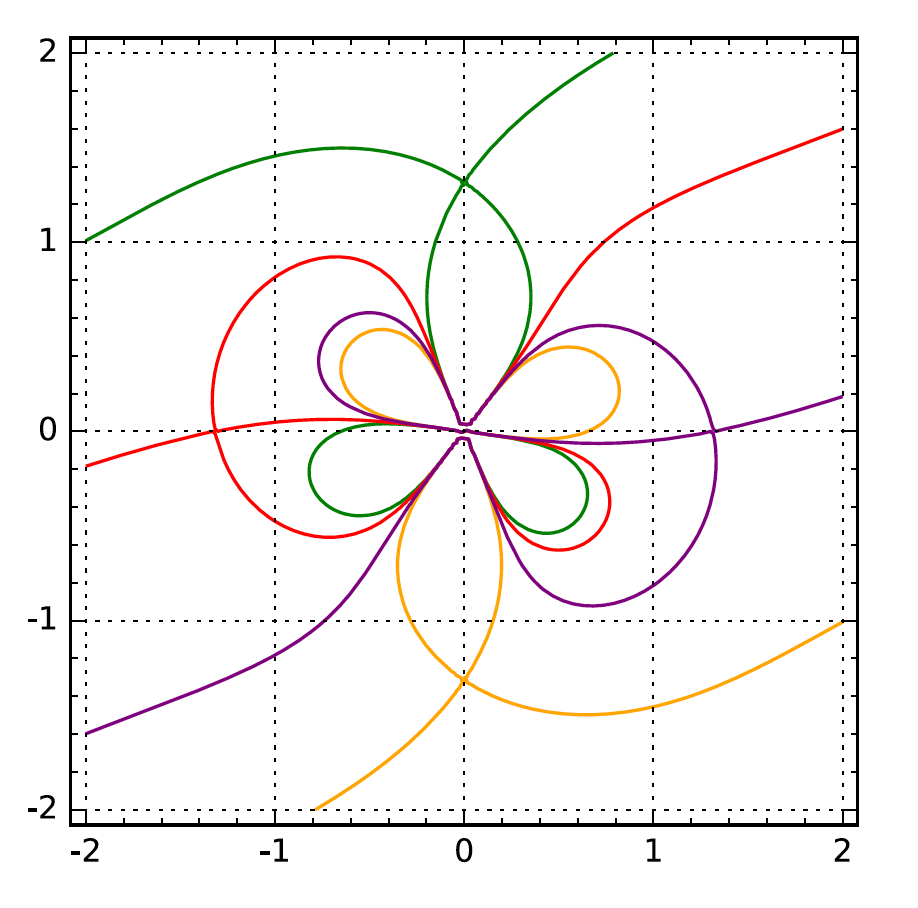}
	\caption{Preimages under $f$}
	\label{preimageslinesP13}
\end{figure}

\begin{figure}
	\hspace*{-9mm}\includegraphics[width=14.5cm]{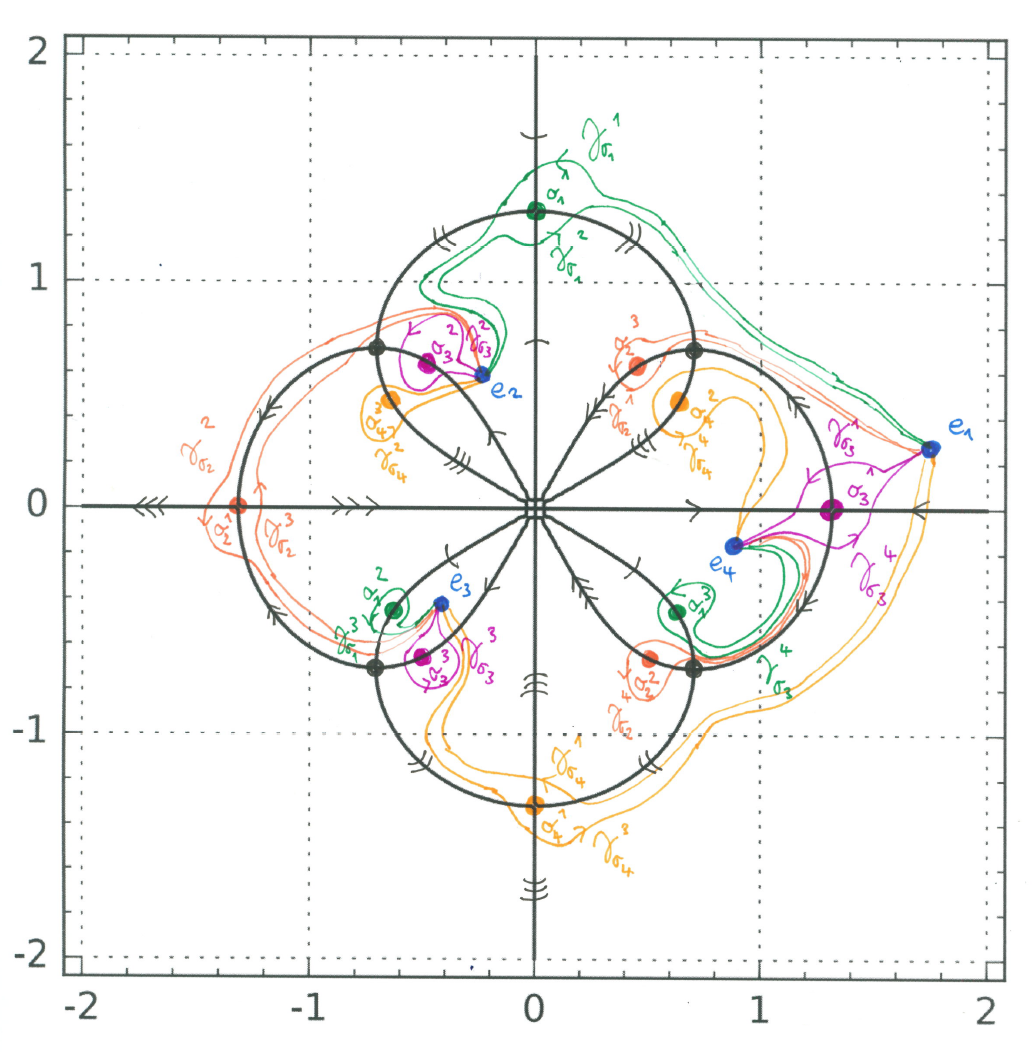}
	\vspace*{-4mm}
	\hspace*{15mm}$\downarrow$
	\newline
	\vspace*{3mm}
	\hspace*{10mm}\includegraphics[width=12cm]{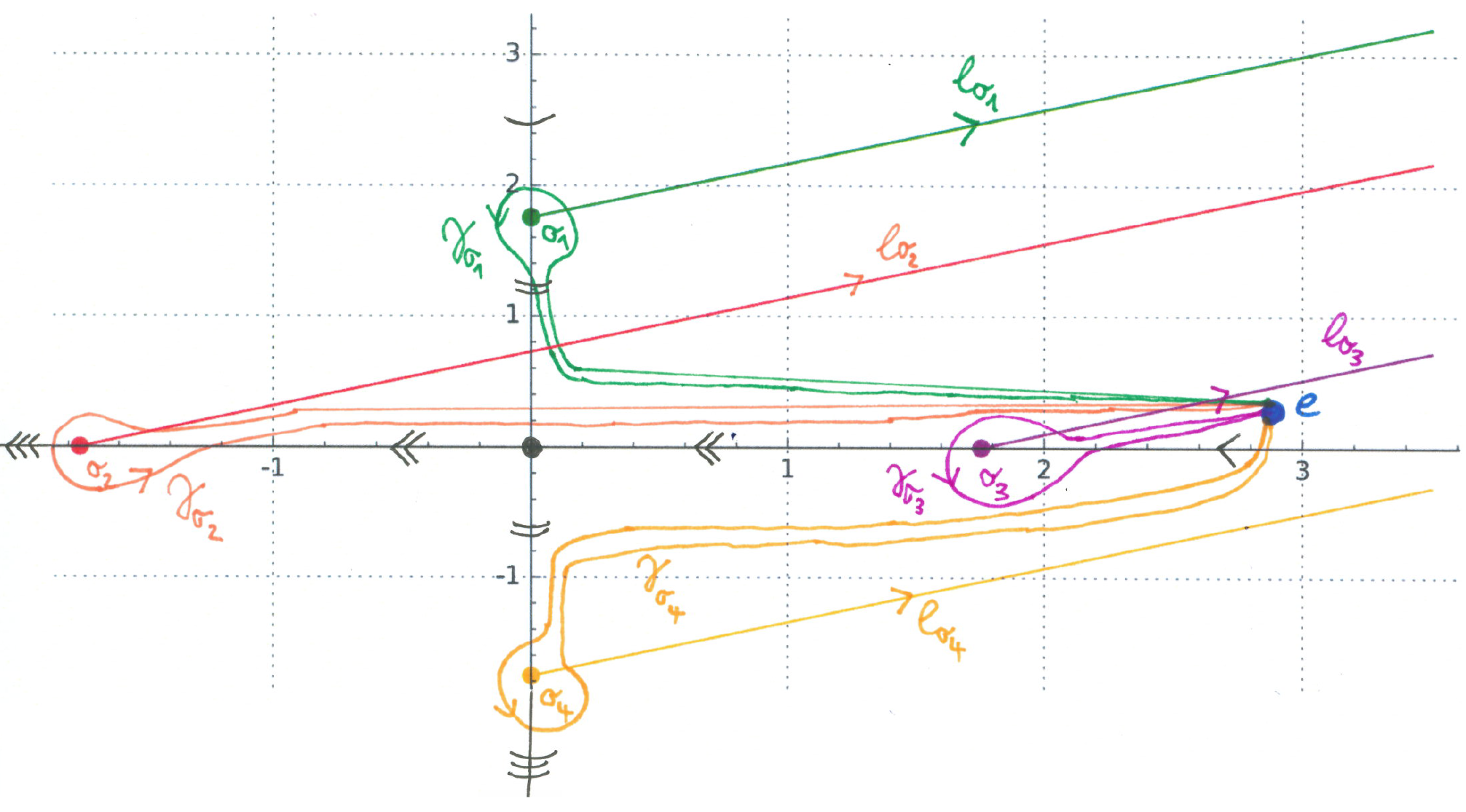}
	\caption{$\gamma_{\sigma_i}$ and their preimages under $f$}
	\label{monodromyP13}
\end{figure}

\begin{figure}
	\hspace*{-9mm}\includegraphics[width=14.5cm]{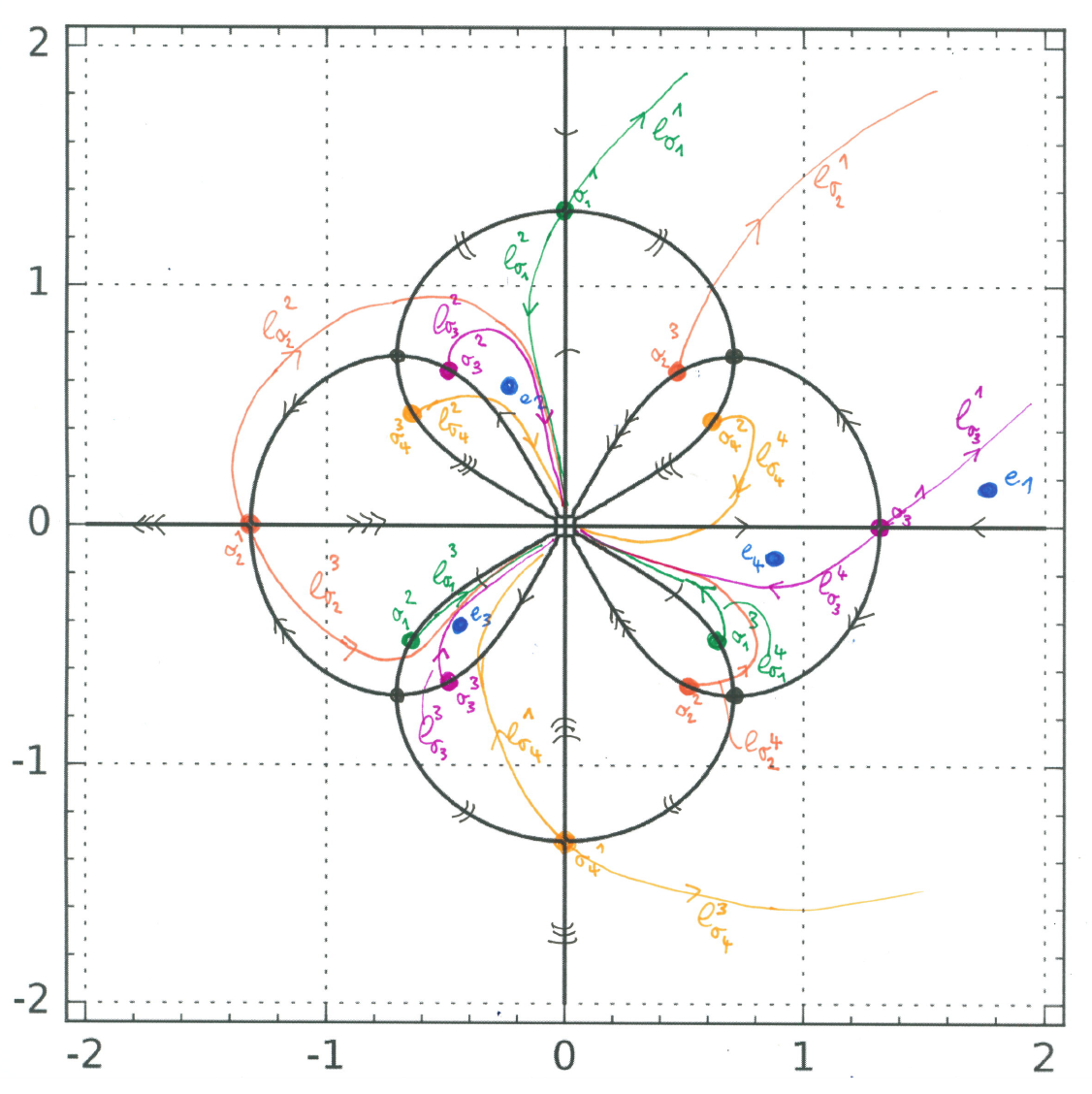}
	\newline
	\hspace*{18mm}$\downarrow$
	\newline
	\hspace*{13mm}\includegraphics[width=12cm]{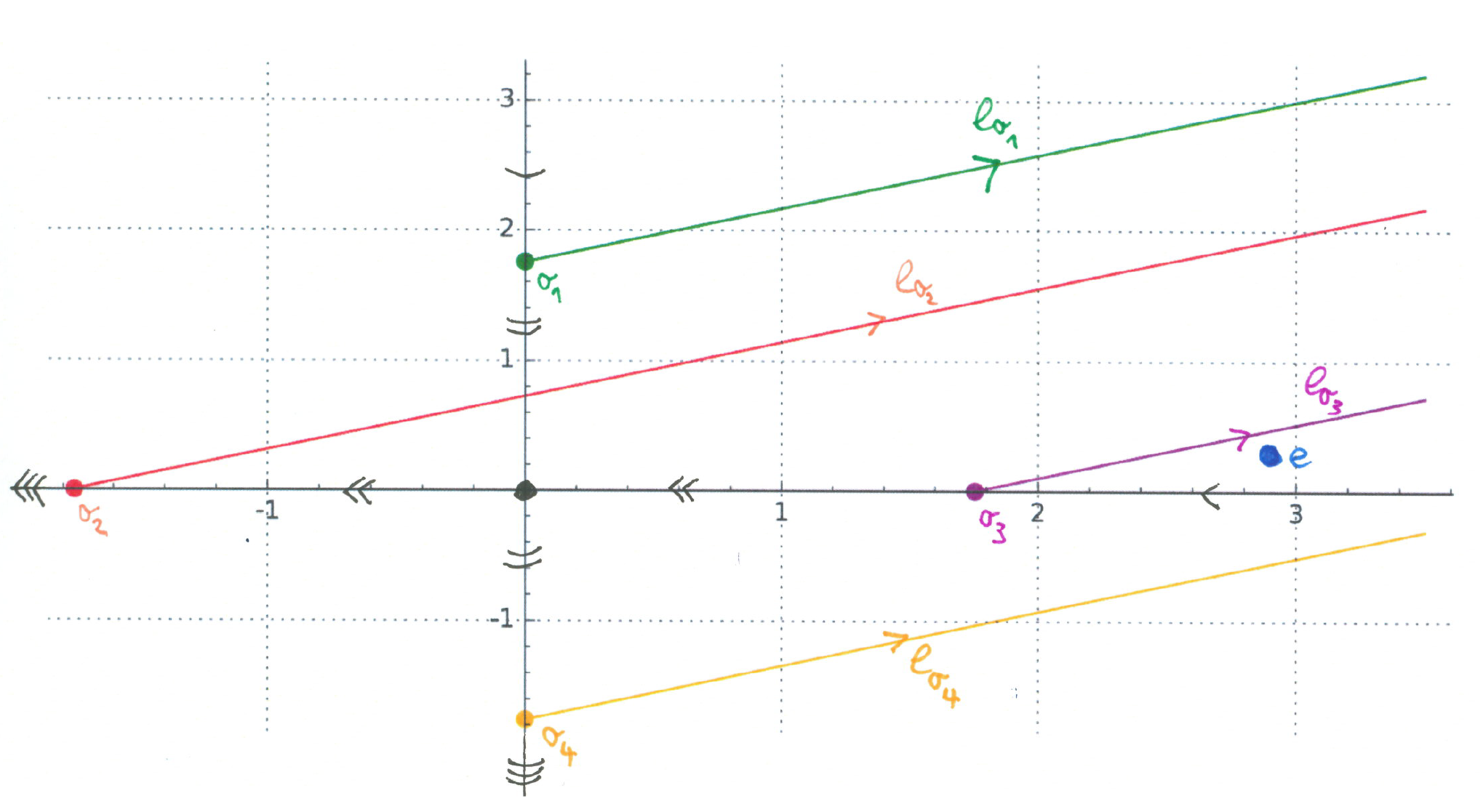}
	\caption{$\ell_{\sigma_i}$ and their preimages under $f$}
	\label{halflinesP13}
\end{figure}

\noindent The blue area in Figure \ref{realimaginarypartyP13} shows where $f$ has real (resp. imaginary) part greater than or equal to $0$. In Figure \ref{preimagesaxesP13}, the preimage of the imaginary (resp. real) axis under $f$ is plotted in blue (resp. red) color. We consider lines passing through the singular values with phase $\frac{\pi }{8}$, as depicted in Figure \ref{linesP13}. The preimages of these lines are plotted in Figure~\ref{preimageslinesP13}. We fix a base point $e$ with $\Re (e) > \Re \left(\sigma_i\right)$ for all $i$ and denote its preimages by $e_1,e_2,e_3,e_4,$ as depicted in Figure \ref{monodromyP13}. In the following, we adopt the notation of\linebreak \cite[Section 4]{DMHS17}. The nearby and global nearby cycles of $F$ are given by~\begin{align*}
\Psi_{\sigma_i}(F) & \coloneqq \R\Gamma_c\left(\A;\mathbb{C}_{\ell_{\sigma_i}^{\times}}\otimes F\right) \simeq H^0\R\Gamma_c \left(\ell_{\sigma_i}^{\times};F\right) \cong \bigoplus_{ e_j \in f^{-1}(e) } \mathbb{C}_{e_j} \cong \mathbb{C}^4,\\
\Psi(F) & \coloneqq \R\Gamma_c\left(\A;\mathbb{C}_{{\A \setminus \ell_{\Sigma}}}\otimes F\right)[1] \simeq \Psi_{\sigma_i}(F) \cong \mathbb{C}^4.
\end{align*}
Furthermore, we fix isomorphisms $i_{\sigma_i}^{-1}F[-1] \cong \bigoplus_{ \sigma_i^j \in f^{-1}\left(\sigma_i\right) } \mathbb{C}_{\sigma_i^j} \cong \mathbb{C}^3.$
\\The exponential components at $\infty$ of the Fourier--Laplace transform of $H^0\left(\int_f \mathcal{O}\right)$ are known to be of linear type with coefficients given by the $\sigma_i\in \Sigma$.
The Stokes rays are therefore given by~$$\left\{ 0,\pm \frac{\pi}{4},\pm \frac{\pi}{2},\pm \frac{3\pi}{4},\pi \right\}.$$
We consider loops $\gamma_{\sigma_i}$, starting at $e$ and running around the singular value $\sigma_i$ in counterclockwise orientation\footnote{counterclockwise orientation since the imaginary part of $\langle \alpha, \beta \rangle$ is positive}, as depicted in Figure \ref{monodromyP13}. We denote by $\gamma_{\sigma_i}^j$ the preimage of $\gamma_{\sigma_i}$ starting at $e_j$, $j=1,2,3,4$.
\\From figure \ref{monodromyP13} we read, in the basis $e_1,e_2,e_3,e_4$, the monodromies 
\begin{align*}
T_{\sigma_1} = \begin{pmatrix}
0 & 1 & 0 & 0\\
1 & 0 & 0 & 0\\
0 & 0 & 1 & 0\\
0 & 0 & 0 &1\\
\end{pmatrix}, \
T_{\sigma_2} = \begin{pmatrix}
1 & 0 & 0 & 0\\
0 & 0 & 1 & 0\\
0 & 1 & 0 & 0\\
0 & 0 & 0 &1\\
\end{pmatrix},\\
T_{\sigma_3} = \begin{pmatrix}
0 & 0 & 0 & 1\\
0 & 1 & 0 & 0\\
0 & 0 & 1 & 0\\
1 & 0 & 0 &0\\
\end{pmatrix},\
T_{\sigma_4} = \begin{pmatrix}
0 & 0 & 1 & 0\\
0 &1 & 0 & 0\\
1& 0 & 0 & 0\\
0 & 0 & 0 &1\\
\end{pmatrix}.
\end{align*}
In order to obtain the maps $b_{\sigma_i}$, we consider the half-lines $\ell _{\sigma_i} \vcentcolon= \sigma_i + \alpha \mathbb{R}_{\geq 0} $. We denote their preimages under $f$ by $\{ \ell_{\sigma_i}^j \}_{j=1,2,3,4}$, depending on which $\gamma_{\sigma_i}^j$ they intersect. We label the preimages of $\sigma_i$ by $\sigma_i^1,\sigma_i^2,\sigma_i^3$, as in Figure \ref{halflinesP13}. The maps $b_{\sigma_i}$ encode which lift of $\ell_{\sigma_i}$ starts at which preimage of $\sigma_i$, induced from the corresponding boundary map in homology.
More explicitly, from Figure \ref{halflinesP13} we read the following:~\begin{itemize}
	\item[\textcolor{green}{$\sigma_1$}:] $\ell_{\sigma_1}^1 \mapsto \sigma_1^1, \ \ell_{\sigma_1}^2 \mapsto \sigma_1^1, \ \ell_{\sigma_1}^3 \mapsto \sigma_1^2, \ \ell_{\sigma_1}^4 \mapsto \sigma_1^3$. \\Therefore, $b_{\sigma_1}$ is the transpose of $\begin{pmatrix} 
	1 & 1 & 0 &0\\
	0 & 0 & 1&0\\
	0 & 0 &0 &1\\
	\end{pmatrix}$.
	\item[\textcolor{red}{$\sigma_2$}:] $\ell_{\sigma_2}^1 \mapsto \sigma_2^3, \ \ell_{\sigma_2}^2 \mapsto \sigma_2^1, \ \ell_{\sigma_2}^3 \mapsto \sigma_2^1, \ \ell_{\sigma_2}^4 \mapsto \sigma_2^2$. \\Therefore, $b_{\sigma_2}$ is the transpose of $\begin{pmatrix}
	0 & 1 & 1 & 0\\
	0 & 0 & 0 &1\\
	1 & 0 & 0 &0\\
	\end{pmatrix}$.
	\item[\textcolor{violet}{$\sigma_3$}:] $\ell_{\sigma_3}^1 \mapsto \sigma_3^1, \ \ell_{\sigma_3}^2 \mapsto \sigma_3^2, \ \ell_{\sigma_3}^3 \mapsto \sigma_3^3, \ \ell_{\sigma_3}^4 \mapsto \sigma_3^1$. \\Therefore, $b_{\sigma_3}$ is the transpose of $\begin{pmatrix}
	1 & 0 &0 &1\\
	0 & 1 & 0&0\\
	0 & 0 &1&0\\
	\end{pmatrix}$.
	\item[\textcolor{orange}{$\sigma_4$}:] $\ell_{\sigma_4}^1 \mapsto \sigma_4^1, \ \ell_{\sigma_4}^2 \mapsto \sigma_4^3, \ \ell_{\sigma_4}^3 \mapsto \sigma_4^1, \ \ell_{\sigma_4}^4 \mapsto \sigma_4^2$. \\Therefore, $b_{\sigma_4}$ is the transpose of $\begin{pmatrix}
	1 & 0 & 1 & 0\\
	0 & 0 & 0 & 1\\
	0 & 1 & 0 & 0\\
	\end{pmatrix}$.
\end{itemize}
We obtain, in the ordered bases $\sigma_i^1,\sigma_i^2,\sigma_i^3$ and $\ell_{\sigma_i}^1,\ell_{\sigma_i}^2,\ell_{\sigma_i}^3,\ell_{\sigma_i}^4$ each:
\begin{align*}
b_{\sigma_1} = \begin{pmatrix}
1 & 0 &0\\
1 & 0 & 0 \\
0 & 1 & 0 \\
0 & 0 & 1\\
\end{pmatrix}, \
b_{\sigma_2} = \begin{pmatrix}
0 & 0 & 1 \\
1 & 0 & 0\\
1 & 0 & 0\\
0 & 1 &0\\
\end{pmatrix},\\
b_{\sigma_3} = \begin{pmatrix}
1 & 0 & 0\\
0& 1 & 0\\
0 & 0& 1\\
1 & 0 &0\\
\end{pmatrix},\
b_{\sigma_4} = \begin{pmatrix}
1& 0 & 0 \\
0 &0 & 1 \\
1 & 0 & 0\\
0& 1 & 0\\
\end{pmatrix}.
\end{align*}
Denote by $u_i\coloneqq u_{\sigma_i}$, $v_i\coloneqq v_{\sigma_i}$, $T_i\coloneqq T_{\sigma_i}$ and $\Phi_i \coloneqq \Phi_{\sigma_i}$. We obtain
\begin{tikzcd}
\!\Phi _i (F) \arrow[r,"v_i"] & \Psi (F) \arrow[l, shift left =.7ex, "u_i"]
\end{tikzcd} 
as the cokernels of the diagrams
\begin{center}
		\begin{tikzcd}
		i_{\sigma_i}^{-1}F[-1] \arrow[r,"b_{\sigma_i}"] \arrow[d]
		& \Psi(F) \arrow[d,"1-T_i"]\\
		0 \arrow[u, shift left =.7ex] \arrow[r]
		& \Psi(F) \arrow[u, shift left =.7ex, "1"]\\	
		\end{tikzcd}
\end{center}
\vspace*{-3mm}
We identify the cokernels of $b_{\sigma_i}$ in the following way:
\begin{itemize}
	\item $\coker b_{\sigma_1} \simeq \mathbb{C} $ via $\left[ \begin{pmatrix}
	v_1\\v_2\\v_3\\v_4
	\end{pmatrix} \right] = \left[ \begin{pmatrix}
	v_1-v_2\\0\\0\\0
	\end{pmatrix} \right],$
	\item $\coker b_{\sigma_2} \simeq \mathbb{C} $ via $\left[ \begin{pmatrix}
	v_1\\v_2\\v_3\\v_4
	\end{pmatrix} \right] = \left[ \begin{pmatrix}
	0\\v_2-v_3\\0\\0
	\end{pmatrix} \right],$
	\item $\coker b_{\sigma_3} \simeq \mathbb{C} $ via $\left[ \begin{pmatrix}
	v_1\\v_2\\v_3\\v_4
	\end{pmatrix} \right] = \left[ \begin{pmatrix}
	v_1-v_4\\0\\0\\0
	\end{pmatrix} \right],$
	\item $\coker b_{\sigma_4} \simeq \mathbb{C} $ via $\left[ \begin{pmatrix}
	v_1\\v_2\\v_3\\v_4
	\end{pmatrix} \right] = \left[ \begin{pmatrix}
	v_1-v_3\\0\\0\\0
	\end{pmatrix} \right].$
\end{itemize}
We obtain that 	
\begin{tikzcd}
(\Phi _i (F) \arrow[r,"v_i"] & \Psi (F) \arrow[l, shift left =.7ex, "u_i"])
\end{tikzcd}
$\simeq $
\begin{tikzcd}
(\mathbb{C} \arrow[r,"v_i"] &\mathbb{C} ^4 \arrow[l, shift left =.7ex, "u_i"]),
\end{tikzcd}
where
\begin{align*}
u_1= \begin{pmatrix}
1 & -1 & 0 &0\\
\end{pmatrix}, \
u_2=\begin{pmatrix}
0 & 1 & -1 & 0\\
\end{pmatrix}, \\
u_3=\begin{pmatrix}
1 & 0 & 0 &-1\\
\end{pmatrix}, \
u_4=\begin{pmatrix}
1 & 0 & -1 & 0\\
\end{pmatrix},
\end{align*} and $v_i=u_i^{\mathrm{t}}$.
Remembering carefully all the choices, by \cite[Theorem 5.2.2]{DMHS17}, we obtain the following Stokes multipliers of the Fourier--Laplace transform of $H^0\left(\int_f \mathcal{O}\right)$ at $\infty$:
\begin{align*}
S_{\beta} = \begin{pmatrix}
1 & u_1v_2&u_1v_3&u_1v_4\\ 0 & 1&u_2v_3&u_2v_4\\
0 & 0 & 1 & u_3v_4\\ 0 & 0 & 0 & 1\\
\end{pmatrix} = \begin{pmatrix}
1 & -1& 1 & 1\\
0 & 1 & 0 & 1\\
0 & 0 & 1 & 1\\
0 & 0 & 0 & 1
\end{pmatrix}, 
\end{align*}~\begin{align*}
S_{-\beta} = \begin{pmatrix}
\mathbb{T}_{1} & 0 & 0 & 0 \\
-u_2v_1 & \mathbb{T}_2 & 0 & 0\\
-u_3v_1 & -u_3v_2 & \mathbb{T}_3 & 0\\
-u_4v_1 & -u_4v_2 & -u_4v_3 & \mathbb{T}_4\\
\end{pmatrix} = \begin{pmatrix}
-1 & 0 & 0&0\\
1 & -1 & 0 & 0\\
-1 & 0 & -1 & 0\\
-1 & -1 & -1 & -1
\end{pmatrix} = -S_{\beta}^{\mathrm{t}},
\end{align*}
where $\mathbb{T}_i \coloneqq 1-u_iv_i$. $S_{\pm \beta}$ describes crossing $h_{\pm \beta}$ from $H_{\alpha}$ to $H_{-\alpha}$, where 
\begin{align*} 
H_{\alpha}= \left\{w \in \left(\A\right)^{\vee} \mid \arg(w) \in\left[- \frac{5 \pi}{8}, \frac{ 3\pi}{8} \right] \right\},\\
H_{-\alpha} = \left\{ w \in \left(\A\right)^{\vee} \mid \arg(w) \in \left[ \frac{ 3 \pi}{8}, \frac{11 \pi}{8} \right] \right\},
\end{align*}
denote the closed sectors at $\infty$ and $h_{\pm \beta}= \pm \mathbb{R}_{>0} \beta \subset \left(\A\right)^{\vee}$, such that\linebreak $ H_{\alpha}\cap H_{-\alpha} = h_{\beta} \cup h_{-\beta}$.

\section{Quantum connection and Dubrovin's conjecture}\label{Gram}
\subsection{Quantum connection}
The quantum connection of a Fano variety (resp. an orbifold) $X$ is a connection on the trivial vector bundle over $\mathbb{P}^1$ with fiber $H^*(X, \mathbb{C})$ (resp. $H_{\mathrm{orb}}^*(X, \mathbb{C})$), the standard inhomogeneous coordinate on $\mathbb{P}^1$ being denoted by~$z$. By \cite[(2.2.1)]{GGI} the quantum connection is the connection given by
\begin{align*}
\nabla_{z \partial_z} = z \frac{\partial }{\partial z} - \frac{1}{z} \left(-K_X \circ \quad\right) + \mu,
\end{align*}
where the first term on the right hand side is ordinary differentiation, the second one is pointwise quantum multiplication by $(-K_X)$, and the third one is the grading operator $$\mu(a) \coloneqq \left( \frac{i}{2} - \frac{\dim X}{2} \right) a \quad \mathrm{for}\ a \in H^i(X, \mathbb{C}).$$ The quantum connection is regular singular at $z=\infty$ and irregular singular at $z=0$.
\\For the weighted projective line $\mathbb{P}(a,b)$, the orbifold cohomology ring is given by (cf. \cite[Example~3.20]{Ma06}) 
\begin{align*} 
H^{\ast}_{\mathrm{orb}}(\mathbb{P}(a,b),\mathbb{C} ) = \mathbb{C} \left[x,y,\xi \right] / \langle xy, ax^{\frac{a}{d}}-by^{\frac{b}{d}}\xi^{n-m},\xi ^d-1 \rangle ,
\end{align*}
where $d=\gcd(a,b)$ and $m,n\in \mathbb{Z}$ such that $am+bn=d$. The grading is given as follows (cf. \cite[Section 9]{AGV}): $\deg x = \frac{1}{A}$, $\deg y = \frac{1}{B}$, $\deg \xi = 0$, where $A=\frac{a}{d} $, $B=\frac{b}{d} $. Quantum multiplication\footnote{We always consider the case $q=1$.} is computed in $$ QH_{\mathrm{orb}}^{\ast} = \mathbb{C} \left[x,y,\xi \right] / \langle xy-1, ax^{\frac{a}{d}}-by^{\frac{b}{d}}\xi^{n-m},\xi ^d-1 \rangle .$$ For $\gcd(a,b)=1$, $-K_{\mathbb{P}(a,b)}$ is given by the element $[x^a+y^b]\in H^1_{\mathrm{orb}}(\mathbb{P}(a,b),\mathbb{C})$. Taking into account that the grading is scaled by $2$, the grading operator is defined by $\mu (a)= \left(i - \frac{\dim X}{2}\right)a$ for $a\in H^i_{\mathrm{orb}}(X,\mathbb{C})$.
\\We obtain the quantum connection of $\mathbb{P}(1,3)$ as follows.
\begin{align*}
H^{\ast}_{\mathrm{orb}}(\mathbb{P}(1,3),\mathbb{C} ) = \mathbb{C} \left[x,y\right] / \langle xy,x-3y^3 \rangle
\end{align*}
with grading given by $\deg x=1,\ \deg y=\frac{1}{3}$. A basis over $\mathbb{C} $ is given by $1,y,y^2,y^3$. Quantum multiplication by $-K_{\mathbb{P}(1,3)}=\left[x+y^3\right]=\left[4y^3\right]$ in this basis is given by the matrix $$ \begin{pmatrix}
0 & \frac{4}{3} & 0&0\\
0 & 0 & \frac{4}{3} &0\\
0 & 0 &0 &\frac{4}{3}\\
4 & 0 & 0 & 0
\end{pmatrix} .$$
The grading $\mu$ is given by the matrix $$\begin{pmatrix}
-\frac{1}{2} & 0 &0 & 0\\
0& -\frac{1}{6} &0 & 0\\
0&0 & \frac{1}{6} & 0\\
0 & 0 & 0 & \frac{1}{2}
\end{pmatrix}.$$
Therefore, the quantum connection of $\mathbb{P}(1,3)$ is given by 
\begin{align}\label{quantum}
\nabla_{z \partial_z} = z \partial_z - \frac{1}{z} \begin{pmatrix}
0 & \frac{4}{3} & 0&0\\
0 & 0 & \frac{4}{3} &0\\
0 & 0 &0 &\frac{4}{3}\\
4 & 0 & 0 & 0
\end{pmatrix} + \begin{pmatrix}
-\frac{1}{2} & 0 &0 & 0\\
0& -\frac{1}{6} &0 & 0\\
0&0 & \frac{1}{6} & 0\\
0 & 0 & 0 & \frac{1}{2}
\end{pmatrix}.
\end{align} It is irregular singular at $z=0$ and regular singular at $z=\infty$. Rewriting in $z^{-1}$ yields the irregular singularity at $\infty$.

\begin{obs}
By the gauge transformation $h=diag\left(\theta^{-\frac{1}{2}},\theta^{-\frac{1}{2}},\theta^{-\frac{1}{2}},\theta^{-\frac{1}{2}}\right)$, which subtracts $\frac{1}{2}$ on the diagonal entries, and passing to $-\theta$, connection \eqref{GM} arising from the Landau--Ginzburg model is exactly the quantum connection \eqref{quantum} of $\mathbb{P}(1,3)$, as predicted by mirror symmetry.
\end{obs}

\subsection{Dubrovin's conjecture}
Let $X$ be a Fano variety (or an orbifold), such that the bounded derived category $\Db(\Coh(X))$ of coherent sheaves on $X$ admits a full exceptional collection $\langle E_1,\ldots, E_n \rangle$, where the collection $\langle E_1,\ldots, E_n\rangle$ is called 
\begin{itemize}
	\item \textit{exceptional} if $\RHom(E_i,E_i)=\mathbb{C} $ for all $i$ and $\RHom\left(E_i,E_j\right)=0$ for $i>j$,
	\item \textit{full} if $\Db(\Coh(X))$ is the smallest full triangulated subcategory of\linebreak $\Db(\Coh(X))$ containing $E_1,\ldots,E_n$.
\end{itemize} In \cite{Dubrovin}, B.~Dubrovin conjectured that, under appropriate choices, the Stokes matrix of the quantum connection of $X$ equals the Gram matrix of the Euler--Poincar\'{e} pairing with respect to some full exceptional collection---modulo some action of the braid group, sign changes and permutations (cf. \cite[Section 2.3]{Cotti}). Then the second Stokes matrix is the transpose of the first one. The Euler--Poincar\'{e} pairing is given by the bilinear form 
\begin{align*} 
\chi (E,F) \coloneqq \sum_k (-1)^k \dim _{\mathbb{C}} \mathrm{Ext}^k(E,F), \quad E,F \in \Db(\Coh(X)).
\end{align*}
The Gram matrix of $\chi$ with respect to a full exceptional collection is upper triangular with ones on the diagonal.
\\For $\mathbb{P}(a,b)$, $\langle \mathcal{O},\mathcal{O}(1), \ldots, \mathcal{O}(a+b-1)\rangle$ is a full exceptional collection of\linebreak $\Db(\Coh(\mathbb{P}(a,b)))$ (cf. \cite[Theorem~2.12]{AKO}). Following~\cite[Theorem~4.1]{CF}, the cohomology of the twisting sheaves for $k\in \mathbb{Z}$ is given by 
\begin{itemize}
	\item $H^0\left( \mathbb{P}(a,b),\mathcal{O}(k)\right) = \bigoplus_{(m,n)\in I_0} \mathbb{C} x^my^n, $ where $$I_0 = \left\{ (m,n) \in \mathbb{Z}_{\geq 0} \times \mathbb{Z}_{\geq 0} \mid am+bn=k \right\}. $$
	\item $H^1 \left( \mathbb{P}(a,b),\mathcal{O}(k) \right)= \bigoplus_{(m,n) \in I_1} \mathbb{C} x^my^n, $ where $$I_1 = \left\{(m,n)\in \mathbb{Z}_{<0} \times \mathbb{Z}_{<0} \mid am+bn=k \right\}. $$
	\item $H^i \left( \mathbb{P}(a,b),\mathcal{O}(k) \right)= 0 $ for all $i\geq 2$.
\end{itemize}
We only need to compute $\mathrm{Ext}^k(\mathcal{O}(i),\mathcal{O}(j))$ for $i<j$ which is given by $H^k\left(\mathcal{O}\left(j-i\right)\right)$ (cf. \cite[Lemma~4.5]{Meier}). Therefore, the zeroth cohomologies of the twisting sheaves $\mathcal{O}\left(j-i\right)$ are the only ones that contribute to the Gram matrix of $\chi$. For $\mathbb{P}(1,3)$ we obtain the cohomology groups $$H^0(\mathcal{O}(1))\cong \mathbb{C},\ H^0(\mathcal{O}(2))\cong \mathbb{C},\ \ H^0(\mathcal{O}(3))\cong\mathbb{C}^2$$ and therefore the Gram matrix of the Euler--Poincar\'{e} pairing with respect to the full exceptional collection $$\mathcal{E}\vcentcolon = \langle \mathcal{O},\mathcal{O}(1),\mathcal{O}(2),\mathcal{O}(3)\rangle$$ is given by
\begin{align}\label{Gram}
S_{\text{Gram}}= \begin{pmatrix}
	1 & 1 & 1 & 2\\
	0 & 1 & 1 & 1\\
	0 & 0 & 1 & 1\\
	0 & 0 & 0 & 1
	\end{pmatrix}.
\end{align}

\subsection{Comparison of the Gram and Stokes matrix}
Mirror symmetry relates the Laurent polynomial $f=x + x^{-3}$ to the weighted projective line $\mathbb{P}(1,3)$. The pair $(\mathbb{G}_m, f=x + x^{-3})$ is a Landau--Ginzburg model of the weighted projective line $\mathbb{P}(1,3)$. According to Dubrovin's conjecture, the Stokes matrix of the quantum connection of $\mathbb{P}(1,3)$ is given by the Gram matrix of the Euler--Poincar\'{e} pairing with respect to some full exceptional collection of $\Db(\Coh(\mathbb{P}(1,3)))$. Note that there is a natural action of the braid group on the Stokes matrix reflecting variations in the choices involved to determine the Stokes matrix (cf. \cite{Gu99}). In our case we have to consider the braid group $$ B_4 = \langle \beta_1, \beta_2, \beta_3 \mid \beta_1\beta_3\beta_1=\beta_3\beta_1\beta_3,\ \beta_1\beta_2\beta_1=\beta_2\beta_1\beta_2,\ \beta_2\beta_3\beta_2=\beta_3\beta_2\beta_3 \rangle .$$
We computed that the Gram matrix of $\chi$ with respect to the full exceptional collection $\mathcal{E} $ is given by \eqref{Gram}. Via the action of the braid $\beta_{1} \in B_4$, we find that it is equivalent\linebreak to $S_{\beta}$. Following \cite[Section~6]{Gu99}, the braid $\beta_1$ acts on the Gram matrix as $$ S_{\text{Gram}} \mapsto S_{\text{Gram}}^{\beta_1} \vcentcolon = A^{\beta_1}\left(S_{\text{Gram}}\right) \cdot S_{\text{Gram}} \cdot \left(A^{\beta_1}\left(S_{\text{Gram}}\right)\right)^{\mathrm{t}},$$
where $A^{\beta_1}\left(S_{\text{Gram}}\right)$ is given by 
\begin{align*}
A^{\beta_1}(S_{\text{Gram}}) = \begin{pmatrix}
0 & 1 & 0 & 0\\
1 & -1 & 0 & 0\\
0 & 0 & 1 & 0\\
0 & 0 & 0 &1
\end{pmatrix} .
\end{align*}
We obtain that $$S_{\text{Gram}}^{\beta_1} = \begin{pmatrix} 
1 & -1 & 1 & 1\\
0 & 1 & 0 & 1\\
0 & 0 & 1 & 1\\
0 & 0 & 0 & 1
\end{pmatrix} =S_{\beta}.$$

\begin{rmk}
$S_{\text{Gram}}^{\beta_1}=S_{\beta}$ is the Gram matrix of the Euler--Poincar\'{e} pairing with respect to the right mutation $\mathbb{R}_1\mathcal{E}$ of the full exceptional collection $\mathcal{E}$ (cf. \cite[Proposition~13.1]{Cotti}). The action of the braid $\beta_1\in B_4$ should correspond to a counterclockwise rotation\linebreak of $\beta$. Therefore, we could expect to have the braid $\beta_1$ acting on our Stokes matrix.
\end{rmk}

\end{document}